\documentclass{amsart}
\usepackage[margin=1.5in]{geometry}

\usepackage[all]{xy}
\usepackage{tikz-cd}
\usepackage{enumerate}
\usepackage{amsfonts}
\usepackage{amssymb, amscd, amsmath}
\usepackage{graphicx}
\usepackage{mathtools}
\usepackage{xcolor}
\usepackage{float}
\usepackage[neverdecrease]{paralist}
\usepackage{url}

\newcommand{\R}{\mathbb{R}}
\newcommand{\C}{\mathbb{C}}

\newcommand{\pf}{\n{\em Proof.}   }

\newcommand{\n}{\noindent}

\newtheorem{teo}{Theorem}[section]
\newtheorem{cor}[teo]{Corollary}
\newtheorem{prop}[teo]{Proposition}
\newtheorem{lema}[teo]{Lemma}

\theoremstyle{definition}

\newtheorem{rmk}[teo]{Remark}

\title{Generalizing quadratic $\mathbb{R}$-Algebraic sets  in $\mathbb{CP}^{n}$ }

\author{Javier Bracho }
\author{Luis Montejano   }

\begin{document}

\maketitle

\quad \quad \quad \quad \quad {\it We dedicate this paper to the memory of our lifelong friend Jorge Arocha.}
\bigskip

\begin{abstract}
Motivated by our study of the complex Banach conjecture, we characterize a complex ellipsoids $\mathcal E$ as compact subsets of $\mathbb C^n$, with the property that every complex line intersect $\mathcal E$ either in a single point or in the complex affine image of the unit disk. This characterization leads to the main interest of this paper. We study the topological behavior of  compact subsets of 
$\mathbb{CP}^n$  with the property that any complex line that intersects them does either at a single point, at the boundary of a complex disk, or along the entire line. In particular, we are interested in quadratic $\R$-algebraic subsets of $\mathbb{CP}^n$.
\end{abstract}

\section{Introduction}
 The images of a ball under complex linear transformations are complex ellipsoids. So, they are the unit balls of finite dimensional Hilbert spaces over the complex numbers.  The analog of segments in the real line, are discs in the complex line. So, the analog of convex bodies in $\mathbb R^n$ should be what we called ``\emph{bombons}'' in $\mathbb C^n$: bodies whose non-trivial sections with complex lines are disks.  Clearly, a complex ellipsoid is a bombon in this sense.  We first introduced bombons in \cite{ABM1} motivated by our study of the complex Banach conjecture \cite{BM}; their name was proposed by our coauthor Jorge Luis Arocha, and alludes to the roundness of  ``{\it bombonas}'', those small gas tanks that were common in Cuba in the mid-20th century.  A little later, they turned out to be  an unexpected characterization of complex ellipsoids with no analog over the real numbers. Namely, in  \cite{ABMCE} we proved that if every complex line that intersects a convex body in $\mathbb C^n$ does so in a disk or a point, then the convex body is a complex ellipsoid.  This characterization naturally leads to the main question of this paper. What are those closed sets of $\mathbb{CP}^n$ with the property that any complex line that intersects them does so either at a single point, at the boundary of a complex disk, or along the entire line? This gives rise to the new definition of bombon that we study in this paper. The analog over the real numbers was studied in \cite{BC}.

\section{Preliminaries}

Let $\C^n$ denote the $n$-dimensional complex vector space, its points are $n$-tuples of complex numbers. 
For $x=\left( x_{1},...,x_{n}\right) \in \mathbb{C}^{n}$, let %
$\bar{x}=\left( \bar{x}_{1},...,\bar{x}_{n}\right) $, where the
bar denotes the complex conjugate. Also, we denote by $\langle x, y\rangle$ the 
interior product of $x$ and $y$; 
$$\langle x, y\rangle= \sum_{i=1}^n x_i \bar{y}_i,$$
(see  \cite{Horn}, Chapter 5),  in such a way that for every $x \in \mathbb{C}^{n}$, $\|x\|^2=\langle x, x\rangle $, where $\|x\|$  is the norm of $x$.
Moverover, denote by $B^{2n}=\{x\in \mathbb{C}^n\mid \|x\|\leq 1\} $, the unit ball.

Adding points at infinity in the classic way, compactify  $\mathbb{C}^{n}$  to obtain the projective $n$-dimensional complex space $\mathbb{CP}^{n}.$   Therefore, $\mathbb{CP}^{1}$ is the Reimann sphere, because it is the one point compactification of the complex line $\mathbb{C}$; that is, $\mathbb{CP}^{1}$ is topologically the $2$-sphere $\mathbb S^2$, and
the projective isomorphisms of $\mathbb{CP}^{1}$ are the classic Moebius transformations.  

We say that $D\subset \mathbb{CP}^{1}$ is a \emph{disk} if $D$ is the image under a Moebius map of the unit disk $\{ z\in \C\mid \|z\|\leq 1\}\subset \mathbb{C}\subset\mathbb{CP}^{1}$.
Similarly, we say that $\Sigma\subset \mathbb{CP}^{1}$ is a \emph{circle} if $\Sigma$ is the image under a Moebius map of the unit circle  $\{ z\in \C\mid \|z\|=1\}\subset\mathbb{C}\subset \mathbb{CP}^{1}$ (please, do not  mistake with the common use of the term circle as a topological $\mathbb S^1$). Some authors use the term \emph{``generalized circles''} for what we are simply calling ``circles''. Thus, the circles in  $\mathbb{CP}^{1}=\C\cup\{\infty\}$ are the ``round'' circles in $\C$ and the closure (in $\mathbb{CP}^{1}$) of real lines in $\C$.  Likewise, the disks of  $\mathbb{CP}^{1}$ are the ``round'' disks in $\C$, the closure  (in $\mathbb{CP}^{1}$) of their complement and the closure (in $\mathbb{CP}^{1}$) of (real) halfspaces in $\C\subset\mathbb{CP}^{1}$.

We say that a closed subset $X\subset \mathbb{CP}^{n}$ is a \emph{bombon} if $\mathbb{CP}^{n}\setminus X$ has exactly two components (or equivalently,  $H^{2n-1}(X,\mathbb Z) = \mathbb Z$) which we call $U$ and $V$, and 
for every complex projective line (which we will simply call \emph{line}) $L\subset  \mathbb{CP}^{n}$: 
$$L\cap X \mbox{ is either the empty set, a single point, a circle or } L\,;$$
and furthermore, if $L\cap X$ is a circle, then one of the components of $L\setminus \Sigma$ is contained in $U$ and the other is contained in $V$. 

\begin{rmk}\label{hypotesis-componentes}
We have decided to include the last statement as a hypothesis because a proof  of this fact, with the ideas of \cite{M}, using the other hypothesis would be intricate and it would divert us from the direction in which we want to conduct this work.
\end{rmk}

Clearly, if $X\subset \mathbb{CP}^{n}$ is a bombon and $\Phi: \mathbb{CP}^{n}\to  \mathbb{CP}^{n}$ is a projective isomorphism, then $\Phi(X)$ is a bombon. Moreover, as we shall prove next, a section of a bombon is either a bombon or a 
\emph{subspace}, that is, a complex projective subspace;  we call it a \emph{$k$-subspace} 
if it has dimension $k$ and a \emph{hyperplane} if furthermore $k=n-1$. 

\begin{lema}\label{sec}
Let $X$ be a bombon in  $\mathbb{CP}^{n}$ and let $H\subset \mathbb{CP}^{n}$ be a subspace.  Then, either $H\cap X$ is a bombon in $H$ or $H\cap X$ is a subspace.
\end{lema}
\pf  By hypothesis, every line $L\subset H$ intersects $H\cap X$ in the empty set, a single point, a circle or in $L$.  If there is a line $L\subset H$ that intersects $X$ in a circle then one of the components of $L\setminus \Sigma$ is contained in $U$ and the other is contained in $V$,
 therefore,  $H\setminus X$ has exactly two components and hence 
$H\cap X$ is a bombon in $H$. If every complex line $L\subset H$ intersects $X$ in the empty set, a single point or in $L$, then,  given two different points $x,y$ of $H\cap X$, the line through $x$ and $y$ is contained in $H\cap X$, and hence $H\cap X$ is a subspace. \qed

\subsection{Examples}\label{examples1}

Now, we briefly introduce two specific examples of bombons and two general ones. Some required proofs which are more technical will be given in the next section.

\smallskip
a) \label{elliptical}If $\mathcal E\subset \mathbb{C}^n \subset  \mathbb{CP}^{n}$ is an $n$-dimensional complex ellipsoid, let $X$ be its boundary. 
Any image of $X$ under a projective isomorphism will be called an \emph{elliptical bombon}. We argued in the introduction that they are bombons. Observe that there is a hyperplane that does not intersect $X$.

\smallskip
b) The closure in  $\mathbb{CP}^{n}$  of a real $(2n-1)$-plane $\Gamma$ 
contained in $\mathbb{R}^{2n}=\mathbb C^n\subset \mathbb{CP}^{n}$ is called the \emph{flat bombon} and is denoted by $\tilde\Gamma$. 

To prove that it is a bombon, first observe that the closure of $
\Gamma$ consists of $\Gamma$ plus all the points at infintity of $\mathbb{CP}^{n}$.  This is so, 
because if $x$ is a point at infinity, it is the point at infinity of an affine complex line
$\ell\subset  \mathbb {C}^n$.  
Either $\ell$ is contained in $\Gamma$ which implies $x\in\tilde \Gamma$, or $\ell$ is parallel to $\Gamma$ and therefore there is $\ell'\subset \Gamma$  parallel to $\ell$ and so $x\in\tilde \Gamma$ or, finally, $\ell\cap \Gamma$ is a real line whose closure in  $\mathbb{CP}^{n}$ contains $x$.  

Consequently, $\mathbb{CP}^{n}\setminus \tilde\Gamma=\mathbb C^n\setminus \Gamma$ which has two contractible components $U$ and $V$
because $\Gamma$ is a real hyperplane of $\mathbb C^n=\mathbb R^{2n}$. Furthermore, if $L$ is a line in $\mathbb{CP}^{n}$, $L$ is either contained in the hyperplane at infinity or $L\cap \mathbb C^n$ is a $2$-dimensional real plane in  $\mathbb C^n=\mathbb R^{2n}$.
Thus, either $L\subset\tilde \Gamma$ or $L$ is parallel to $\Gamma$ and therefore $L\cap\tilde\Gamma$ is a single point at infinity, or $L\cap \Gamma$ is a real line $R$.  This last case implies that $L\cap \tilde\Gamma$ is a circle and  one of the components of $L\setminus R$ is contained in $U$ and the other is contained in $V$. 
All this proves that the flat bombon $\tilde\Gamma$ is indeed a bombon.  Note that a line never intersects the flat bombon in the empty set.

\smallskip
c) We first need a general definition. Given two non-void disjoint subsets $A$ and $B$ of $\mathbb{CP}^{n}$, denote by $A\star B$ the union of all the lines in $\mathbb{CP}^{n}$ that \emph{intersect} (that is, that have non-empty intersection with) both $A$ and $B$. We will call $A\star B$ the \emph{join} of $A$ and $B$;  beware that it is not the topological ``join'' with segments, it is related but it is a union of complete complex projective lines parametrizaed by the cartesian product. 

Let $X$ be a bombon in $\mathbb{CP}^{n-k-1}$ and let $\Delta\subset \mathbb{CP}^{n}$ be a 
$k$-subspace of $\mathbb{CP}^{n}$ (recall that $\Delta$ is projectively $\mathbb{CP}^{k}$) 
such that $\Delta\cap \mathbb{CP}^{n-k-1}=\emptyset$.  Then, $X\star \Delta$ is a bombon called the \emph{join of $X$ with apex $\Delta$}. 

Clearly the  complement of $X\star \Delta$ in $\mathbb{CP}^{n}$ has two connected components because
the complement of $X$ in $\mathbb{CP}^{n-k-1}$ has two components.  If $L$ is a line that intersects $\Delta$, it is easy to see that it intersects $X\star \Delta$ in either a point or in $L$. If not, then its projection from $\Delta$ to $\mathbb{CP}^{n-k-1}$ is a line $L'\subset\mathbb{CP}^{n-k-1}$ that intersects $X$ as  $L$ does $X\star \Delta$.

Finally, we will say that $X\subset \mathbb{CP}^{n}$  is a  \emph{conical bombon} if $X$ is the join of an elliptical bombon with a 
$k$-subspace, $0\leq k\leq n-1$. 

Note that in $\mathbb{CP}^n$, two elliptical bombons are projectively equivalent. This immediately implies that  two conical bombons  which are the join of an elliptical bombon with a complex $k$-subspace, $k<n-2$, are projectively equivalent. 

\smallskip
d) Consider homogeneous coordinates in $\mathbb{CP}^{n}$  and let $X\subset \mathbb{CP}^{n}$ be the following set 
\begin{equation}\label{BombonaAlgebraica}
X=\{[x_0:x_1:\dots: x_n]\in\mathbb{CP}^n \,\mid\,  \sum_{j=0}^n \epsilon_j\|x_j\|^2=0\}
\end{equation}
where $\epsilon_j\in\mathbb{R}$ and at least one \emph{coefficient} $\epsilon_j$ is positive and one negative.   
Any projective image of $X$  will be called an \emph{algebraic bombon}, and is what we referred to as a quadratic $\R$-algebraic set in the title. 

Let $X$ be as in \eqref{BombonaAlgebraica}. In the next section we prove that it satisfies the intersection propperties with lines of a bombon. Meanwhile, observe that $\mathbb{CP}^n\setminus X$ has exactly two components defined by whether the defining expression of $X$ is positive or negative, because of our assumption on the coefficients $\epsilon_i$.  Moreover, we have nonempty index subsets $I$ and $J$ of cardinalities $p+1$ and $q+1$ for which the respective coeficients are positive or negative, and we may assume that $0\leq p\leq q$. We will call $X$ an \emph{algebraic bombon} of type $(p,q)_n$. 

Note that elliptical bombons are algebraic bombons of type $(0,n-1)_n$, flat bombons are algebraic of type $(0,0)_n$ and conical bombons are algebraic of type $(0,q)_n$ with $1\leq q <n-1$.

\subsection{Singular locus, cores and full bombons.}

Let $X\subset \mathbb {CP}^n$ be a bombon. We say that $x\in X$ is \emph{a singular point of} $X$ if every complex line $L$ through $x$ intersects $X$ either in the single point $\{x\}$ or in $L$.   The set of all singular points of $X$, denoted $\Sigma(X)$ is its \emph{singular locus}. 

\begin{lema}\label{lema:SingularLocus}
Let $X\subset \mathbb {CP}^n$ be a bombon. 
Then, its singular locus $\Sigma(X)$ is a subspace.
\end{lema}

\pf Consider $x,y\in\Sigma(X)$, $x\neq y$, and let $L$ be the line through them; observe that $L\subset X$. It is enough to prove that $L\subset\Sigma(X)$. Let $z\in L$ be different from $x$ and $y$, and let $L^\prime$ be a complex line through $z$, different from $L$. To see that $z\in\Sigma(X)$, we need to prove that $L^\prime \cap X$ is either $\{z\}$ or $L^\prime$. Suppose that there exists $z^\prime\in L^\prime\cap X$, $z^\prime\neq z$; we must show that $L^\prime \subset X$ to complete the proof. 

Consider the \emph{plane} (i.e., the $2$-subspace) $P$, generated by $L$ and $L^\prime$, and observe that it is the plane generated by $x$, $y$ and $z^\prime$. Let $L_y$ be the line through $y$ and $z^\prime$. Then, because $y$ is a singular point of $X$, $L_y\subset X$ because it contains two points in $X$. Now, because $x$ is also singular, all the lines from $x$ to points in $L_y$ are in $X$. But their union is $P$, so that $L^\prime\subset P\subset X$.  
\qed 

\smallskip

Let $X\subset \mathbb {CP}^n$ be a bombon.  We say that $X$ is a \emph{smooth bombon} if $X$ does not contain any singular point, that is if $\Sigma(X)=\emptyset$.

From Lemma~\ref{lema:SingularLocus}, it easily follows that if 
$X\subset \mathbb {CP}^n$ is a bombon and $\Gamma$ is a complementary  subspace to $\Sigma(X)$, 
then $X'=X\cap \Gamma$ is a smooth bombon and $X= X'\star \Sigma(X)$. That is: 

\begin{lema}\label{SingJoinSmooth}
Let $X\subset\mathbb{CP}^n$ be a bombon. If $X$ has singular points, then $X$ is the join of a smooth bombon $X'$ with apex its singular locus $\Sigma(X)$.     \qed
\end{lema}

Let $X\subset \mathbb {CP}^n$ be a bombon with $U$ and $V$ the components of $\mathbb {CP}^n\setminus X$. We say that a subspace $\mathcal C_U$ (respectivelly, $\mathcal C_V$) is a \emph{core} of $X$ if $\mathcal C_U\subset U$ (resp., $\mathcal C_V\subset V$)
and it has maximum dimension among all such subspaces contained in $U$ (resp., $V$). Let $p$ and $q$ be the dimesions of the cores $\mathcal C_U$ and $\mathcal C_V$ of $X$ (and we may further assume that $p\leq q$), then we say that $X$ is a bombon of \emph{type} $(p,q)_n$. 

Observe that the type of algebraic bombons, as defined previously, corresponds to this general definition.

We say that a bombon $X$ of type $(p,q)_n$ is \emph{full} if 
$$p+q+dim(\Sigma(X))=n-2\,,$$
where $dim(\emptyset)=-1$. That is, if the two cores and the singular locus generate $\mathbb{CP}^n$. 
Consequently, a full smooth bombon has type $(p,q)_n$, with $p+q=n-1$.

\smallskip

The following lemma allows the use of the theory of fiber bundles  in the study of the structure of bombons.  First some definitions.  
Let $\gamma_1^n$ be the canonical complex line bundle, that is:  
$\gamma_1^n= \{(\ell, x)\mid \ell \mbox{ is a line in } \mathbb {C}^n \mbox{ through the origin and } x\in \ell\}$.
Consider also the corresponding unitary $\mathbb S^1$-bundle
$S\gamma_1^n= \{(\ell, x)\mid \ell \mbox{ is a line in } \mathbb {C}^n \mbox{ through the origin and } x\in \ell, \|x\|=1\}$.  

\begin{lema}\label{cfb}
Let $X\subset \mathbb {CP}^n$ be a bombon with $U$ and $V$ the two components of $\mathbb {CP}^n \setminus X$. Let $z\in U$ and let $\mathcal L_z$ be the space of lines through $z$.  Suppose that for every $L\in \mathcal L_z $, $L\cap X$ is a circle. 
If $\mathcal E_z=\{ (L,x)\in \mathcal L_z \times \mathbb {CP}^n\mid x \in L\cap X\}$ and the projection $\pi:\mathcal E_z\to \mathcal L_z$ is given by $\pi(L, x)=L$, then 
$\pi:\mathcal E_z\to \mathcal L_z$ is equivalent to the $\mathbb S^1$-bundle $S\gamma_1^n$.
\end{lema}

\pf First of all note that $\{ (L,x)\in \mathcal L_z \times \mathbb {CP}^n\mid x \in L\}$ is a fiber bundle with fiber $\mathbb{CP}^1$ and structure group the Moebius maps. Consider the natural metric on $\mathbb{CP}^n$, known as the Fubini-Study metric. It endows each line with a classic ``round'' metric of $\mathbb{S}^2$. Let $B\subset U$ be a small ball centered at $z\in \mathbb {CP}^n$ with boundary $S\subset U$, naturally identified with $\mathbb{S}^{2n-1}$.   The corresponding fibers of $L\in\mathcal L_z$ are two circles, one is $L\cap X=\Sigma$ and the other $L \cap S$ centered at $z$. Consider the closure of the two components of $L\setminus\Sigma$, which are disks. One of them is contained in $U$, and has a well defined center $c_L\in L\cap U$. Now, we can specify a Moebius map on $L$ that sends $S\cap L$ to $X\cap L= \Sigma$. There is a well defined ``radius''  (within a geodesic great circle) from $z$ to $c_L$ in $L \cap U$. Translate $z$ to $c_L$ along the radius (which is a rotation arround the corresponding poles) followed by the appropriate real dilation which fixes $c_L$ and sends (the translate of) $S\cap L$ to $\Sigma$. 
Since $c_L$ and the radius vary continuously with $L$, this gives an equivalence between the $\mathbb S^1$-bundles $\mathcal E_z$ and $S\gamma_1^n$. \qed

\section{Examples of Bombons}

In this section we first complete the proof that algebraic bombons are bombons and then characterize flat, elliptical and conical bombons in terms of their behaviour with respect to subspaces.  

\subsection{The algebraic bombons}

\begin{teo}\label{AlgBomb}
Let $X\subset \mathbb{CP}^{n}$ be an algebraic bombon.  Then $X$ is a full bombon. 
\end{teo}

\pf Recall that an algebraic bombon $X$ is defined by equation \eqref{BombonaAlgebraica}.
We first complete the proof that 
$X$ is indeed a bombon. 
We have already argued that  $X$ separates $\mathbb{CP}^{n}$ in two components, according to whether the outcome of the defining expression, $\sum_0^n \epsilon_j\|x_j\|^2$, 
is positive or negative (it is $0$ in $X$), because of the hypothesis that at least one  \emph{coefficient} $\epsilon_j$ is positive and one is negative.

Let $L$ be a line that intersects $X$ in more than a single point. Let $a=[a_0:a_1:\dots: x_n]$ and $b=[b_0:b_1:\dots: b_n]$ be two different points in $L\cap X$. Then, $L$ can be parametrized by $[x:y]\in\mathbb{CP}^1$:
$$[x:y]\mapsto[xa_0+yb_0:\dots:xa_n+yb_n]\,.$$

The defining expression of $X$ yields for these points:
\begin{align*}
 \sum_{j=0}^n \epsilon_j\|xa_j+yb_j\|^2 =& 
 \sum_{j=0}^n \epsilon_j(xa_j+yb_j)\overline{(xa_j+yb_j)} 
 =\sum_{j=0}^n \epsilon_j(xa_j+yb_j)(\bar{x}\bar{a}_j+\bar{y}\bar{b}_j)\\
 =& \sum_{j=0}^n\epsilon_j(a_j\bar{a}_j x\bar{x} + b_j\bar{b}_j y\bar{y} + a_j\bar{b}_j x\bar{y}+b_j\bar{a}_j y\bar{x})\\
 =& (\sum_{j=0}^n\epsilon_ja_j\bar{a}_j)x\bar{x}
 +(\sum_{j=0}^n\epsilon_j b_j\bar{b}_j)y\bar{y}
 +(\sum_{j=0}^n\epsilon_j a_j\bar{b}_j)x\bar{y}
+(\sum_{j=0}^n\epsilon_j \bar{a}_jb_j)\bar{x}y\,.
\end{align*}

The first two summands are $0$ because of our choice of points in $X$. So, to describe $L\cap X$ we have to find the solutions $[x:y]\in\mathbb{CP}^1$ of the equation
\begin{equation}\label{eq1}
cx\bar{y}+ \bar{c}\bar{x}y=0\,,\quad\text{where}\quad
c=\sum_{j=0}^n\epsilon_j a_j\bar{b}_j\,.
\end{equation}
There are two cases. First, if $c=0$ then $L\cap X=L$ because any $[x:y]\in\mathbb{CP}^1$ makes the expression $0$. And second, if $c\neq 0$, observe that the expression $x+\bar{x}$ in $\C$ takes all possible real values and the equation $x+\bar{x}=0$ in $\C$ defines the imaginary axis $\{si\,\mid\,s\in\R\}$. So that we may parametrize the solutions of \eqref{eq1} by the real projective line $[s:t]\in\mathbb{RP}^1$ (with $s,t\in\R$), because taking $[x:y]=[si:tc]$ in \eqref{eq1}, we have
$$c(si)(\overline{tc})+ \bar{c}(\overline{si})tc=csit\bar{c}-\bar{c}sitc=0\,.$$
Therefore, $L\cap X$ is a circle and the two sides of it make the defining expression positive and negative respectivelly. We have proved that any line $L$ in $\mathbb{CP}^n$ intersects $X$ in either the empty set, a single point, a circle or $L$, so $X$ is a bombon.

We are left to prove that $X$ is a full bombon. Let $I$ and $J$ be the sets of indices, $0\leq j\leq n$,  for which $\epsilon_j$ has the same sign, and $K$ be the remaining ones (for which $\epsilon_j=0$). Clearly and without loss of generality, we may assume that $p:=\#I-1\leq\#J-1=:q$, and by hypothesis, we have that $0\leq p, q$. Then, the subspaces generated by the coordinates with indices in $I$ and $J$ (of dimensions $p$ and $q$ respectivelly) are cores of $X$; and the subspace generated by the coordinates in $K$ (that do not appear in the defining expression) generate its singular locus. So that $X$ is a full bombon of type $(p,q)_n$.
\qed 

\medskip
In particular, we have proved the existence of full bombons of type $(p,q)_n$ for any $p,q\geq 0$ such that $p+q\leq n-1$. And moreover, it is easy to produce a  matrix in $\text{Gl}(n+1,\C)$ to prove:

\begin{teo}\label{teoAlgsEquiv}
Two algebraic bombons of the same type are projectively equivalent. \qed
\end{teo}

\medskip
\subsection{The flat bombons}

The purpose of this subsection is to prove that flat bombons are precisely those bombons that contain a complex hyperplane. 

\begin{teo}\label{thmdeg}
Let $X\subset \mathbb{CP}^{n}$ be a bombon. The bombon $X$ is a flat bombon if and only if there is a complex hyperplane contained in $X$
\end{teo}

\pf Let $H$ be a complex hyperplane $H$ of $\mathbb{CP}^{n}$ and suppose $H\subset X$.   Without loss of generality we may assume that $H$ is the hyperplane at infinity, so  
$\mathbb{CP}^{n}\setminus H=\mathbb C^{n}$. Our next purpose is to prove that $\mathbb C^{n}\cap X$ is a $(2n-1)$-dimensional plane. Let $x$ and $y$ two different points in 
$\mathbb C^{n}\cap X$ and let $L$ be the complex line through $x$ and $y$. Therefore, either $L\cap H$ is a circle or $L$. In both cases, the real line in $\mathbb C^{n}=\mathbb R^{2n}$, through $x$ and $y$ is contained in $X$.
This implies that $\mathbb C^{n}\cap X$ is a real affine subspace of $\mathbb C^{n}=\mathbb R^{2n}$.  Since $\mathbb{CP}^{n}\setminus X= \mathbb C^{n}\setminus X$ has two connected components, we conclude that
$\mathbb C^{n}\cap X$ is a real $(2n-1)$-dimensional subspace $\Gamma\subset \mathbb C^{n}$ and consequently that $\tilde\Gamma=X$. \qed

\smallskip

Note that the complement of a flat bombon are two contractible sets.

\begin{cor}\label{cordeg}
Any flat bombon  in $\mathbb{CP}^{n}$ is the algebraic bombon of class $(0,0)_n$. \qed
\end{cor}

And, in view of Theorem~\ref{teoAlgsEquiv}, we also have:

\begin{cor}\label{cordeg}
Any two flat bombons in $\mathbb{CP}^{n}$ are projetively equivalent.  \qed
\end{cor}

\medskip
\subsection{The Elliptic Bombons I}\label{secEllBomI}
Let us define that a bombon $X$ is \emph{elliptic} if it contains no lines. Thus, it is a closed subset $X\subset \mathbb{CP}^{n}$ such that $\mathbb{CP}^{n}\setminus X$ has exactly two components $U$ and $V$ and 
for every line $L\subset  \mathbb{CP}^{n}$, 
$L\cap X$  is either the empty set, a single point or a circle, where  one of the components of $L\setminus \Sigma$ is contained in $U$ and the other is contained in $V$.

The purpose of this section is to prove that  what we called elliptical bombons in \ref{examples1}.(a), that are also the algebraic bombons of type $(0,n-1)_n$, are precisely the elliptic bombons as defined by their intersection with lines.

\begin{teo}\label{thmEB}
Let $X\subset \mathbb{CP}^{n}$ be an elliptic bombon. Then $X$ is the boundary of a complex ellipsoid in $\C^n\subset\mathbb {CP}^n$. 
\end{teo}

The proof is by induction and will be completed until the next section. For the moment, let us make some general remarks and advance on it.

Note that an elliptic bombon in $\mathbb{CP}^{1}$ is a circle. Hence, the theorem is true for $n=1$.  Suppose the theorem is true for $n-1$. Let $X\subset \mathbb{CP}^{n}$ be  an elliptic bombon and suppose $H$ is a complex hyperplane.  Since $X$  does not contain lines, then
by Lemma \ref{sec} $H\cap X$ is either the empty set, a single point or a bombon in $H$. Furthermore, by induction, if $H\cap X$ is a bombon, then it must be the boundary of a complex $(n-1)$-ellipsoid contained in $H$.

Accordingly, the proof of Theorem \ref{thmEB} has three cases:
\begin{enumerate}
\item there is a complex hyperplane $H$ that does not intersect $X$,
\item there is a complex hyperplane $H$ that intersects $X$ in a single point,
\item every complex hyperplane $H$ intersects $X$ in  the boundary of a complex $(n-1)$-ellipsoid contained in $H$.
\end{enumerate}

\noindent\emph{The proof of Case} (1). Suppose there is a complex hyperplane that does not intersect $X$, Then we may assume without loss of generality that  it is the hyperplane at infinity and assume, therefore,
that $X\subset \mathbb C^n$.  Note that, by hypothesis $\mathbb C^n\setminus X$ has two components, one of which is not bounded and will be denoted by $V$ and the other which is bounded and denoted by $U$.

We first prove that $X\cup U$ is a convex set $K$ whose boundary is $X$.

Consider $a \not= b \in K$ and let $L$ be the line through $a$ and $b$. Then, $L\cap X$ is a circle, because it contains two different points.   
Therefore, $L\cap K$ is a disk that contains $a$ and $b$ and consequently the closed segment between $a$ and $b$ is in $K$.  

To complete the proof of Theorem \ref{thmEB}, in the case in which there is a complex hyperplane that does not intersect $X$,  it is enough to use the following theorem whose proof  
 is in \cite{ABMCE}.  For completeness, we include an alternative proof in the Appendix. 

\begin{teo}\label{thmB}
Let $K\subset \mathbb C^n=\mathbb R^{2n}$ be a convex body. Then $K$ is a complex ellipsoid if and only if  every complex line intersects $K$ either in the empty set, a single point or a disk.
\end{teo}

\medskip

\noindent\emph{The proof of Case} (2). Now let us consider the case in which there is a hyperplane that intersects $X$ in a single point.  Suppose that $H$ is a hyperplane with the property that $H\cap X=\{\infty\}$ and assume, without loss of generality, that 
$H$ is the hyperplane at infinity of  $\mathbb {CP}^n$.  Consequently, $X\setminus \{\infty\}\subset \C^n$ is a closed, unbounded subset of $\C^n=\R^{2n}$ with the property that the complement of $X$ in $\mathbb C^n$ has two components.
One of them, called $V$, contains the hyperplane at infinty minus $\{\infty\}$. The other, contained in $\mathbb C^n$ is $U$.   

We will prove that $U$ is an  open convex set in  $\mathbb C^n=\mathbb R^{2n}$. Let $a \not= b \in U$ and let $L$ be the line through $a$ and $b$. Suppose first that $\infty$ is not in $L$, then $L\cap X$ is a circle, because $L$ intersects $V$ at infinity.   
Therefore, $L\cap U$ is an open disk in $L\cap\C^n$ that contains $a$ and $b$ and consequently the closed segment between $a$ and $b$ is also in $U$. 

If $\infty\in L$, then $L\cap X\cap\C^n$ is a real line in $L\setminus\{\infty\}$ with $a$ and $b$ in the same open half-plane. Then, the segment between $a$ and $b$ is in $U$, completing the proof that $U$ is convex. 

Let $K\subset\C^n$ be the closure of $U$ in $\C^n$. It can also be defined as $K:=(X\setminus \{\infty\}) \cup U$, and it is a convex body in $\R^{2n}=\C^n$.
 Note that $\text{int}(K)=U$ and $\text{bd}(K)=X\setminus \{\infty\}$. Now, consider $w\in \mbox{bd}(K)$.  Within $\C^n=\R^{2n}$, let $\Delta$ be a real support hyperplane of $K$ at $w$, and observe that $\Delta\cap U=\emptyset$. Let $\Gamma$ be the unique affine complex hyperplane contained in $\Delta$ and containing $w$. Denote by $\tilde\Gamma$ the hyperplane in $\mathbb {CP}^n$ containing $\Gamma$.  

If $\infty\in\tilde \Gamma$, let $L$ be the line through $w$ and $\infty$. Since $L\setminus\{\infty\}$ is contained in $\Gamma\subset\Delta$, we have that $L\cap U=\emptyset$, and therefore, $L\cap X$ is not a circle. But $L$ has at least two points in $X$, so it is contained in $X$, contrary to the definition of elliptic bombon. 
Therefore, $\infty$ is not a point of $\tilde \Gamma$. Consider a parallel real hyperplane $\Delta^\prime$ to $\Delta$ that does not intersect $K$ (contained in $V$). It contains a complex hyperplane $\Gamma^\prime$, parallel to $\Gamma$, whose closure in $\mathbb{CP}^n$ does not intersect $X$, because it  does not contain $\infty$. Therefore, we are in Case 1.  

To summarize, we have proved the following:

\begin{prop}\label{propEll} 
Let $X$ be an elliptic bombon in $\mathbb{CP}^n$. If there exists a hyperplane $H$ such that $X\cap H$ is either empty or a single point, then $X$ is the boundary of a complex ellipsoid in $\C^n\subset\mathbb{CP}^n$.   \qed
\end{prop}
 
\medskip

We will eventually see that Case 3 is impossible, but for now we can prove it is so in dimension $2$. 

\begin{teo}\label{teoDim2}
Let $X$ be a bombon in $\mathbb{CP}^2$. Then either it is the flat bombon or the boundary of a complex ellipsoid in $\C^2\subset\mathbb{CP}^2$, that is, it is an algebraic
bombon of type $(0,0)_2$ or $(0,1)_2$ respectively.
\end{teo}

\pf
Suppose that  $X\subset \mathbb{CP}^{2}$ is a bombon. If $X$ contains a line, it is the flat bombon by Theorem~\ref{thmdeg}, because lines are hyperplanes in $\mathbb{CP}^2$. Therefore, we are left to assume that $X$ is elliptic. In view of Proposition~\ref{propEll}, we must also assume that $X$ has the property that for every line $L$, $L\cap X$ is a circle (that is, that $X$ falls in Case 3). Let $U$ and $V$ be the components of $\mathbb{CP}^{2}\setminus X$. We shall prove that
$U\cup X$ and $V\cup X$ are contractible spaces to conclude that this cannot happen. 

Let $z\in U$. Since every line $L$ through $z$ intersects $X$ in a circle.  By Lemma \ref{cfb}, 
 the canonical bundle  $\mathcal E_z=\{ (L,x)\in \mathcal L_z \times \mathbb {CP}^2\mid x \in L\cap X\}$  is equivalent to the unitary canonical $\mathbb S^1$-bundle $S\gamma_1^n$.  Consequently, $U\cup X$ is contractible and similarly $V\cup X$ is contractible, 
 thus proving that the Lusternik-Schnirelman category of $\mathbb {CP}^2$ is two, which is a contradiction. \qed

\section {The Tangent Space and Elliptic Bombons II}

To conclude the proof of Theorem~\ref{thmEB}, we need a general definition.

Let $X\subset \mathbb{CP}^n$ be a bombon and $x\in X$.  Define the \emph{tangent space of $X$ at $x$} for $n>1$ as:
$$T_xX=\{ y\in \mathbb{CP}^n\mid y \mbox{ belongs to a line } L \mbox{ through } x, \mbox{ where }  L\cap X \mbox{ is not a circle}\}\,.$$
And to complete the definition, for $n=1$ define $T_xX=\{x\}$.

It is clear that  $x$ is a singular point of $X$ if and only if   $T_xX=\mathbb{CP}^n$. And also, that for every $x\in X$, the singular locus is contained in its tangent space, that is $\Sigma(X)\subset T_xX$. 

\begin{lema}\label{tangent}
Let $X\subset \mathbb{CP}^n$ be a bombon and let $x\in X$ be a non singular point, then $T_xX$ is a hyperplane.
\end{lema}

\pf 
Since the case $n=1$ is trivially true (points are hyperplanes of $\mathbb{CP}^1$), we start by proving the case $n=2$. 

By Theorem~\ref{teoDim2}, we have two cases. If $X$ is the conical bombon of type $(0,0)_2$ and $x\in X$ is not its unique singular point, $z$ say, then the line trhrough $x$ and $z$ (in $X$) is $T_xX$ because any other line through $x$ cuts $X$ in a circle. If $X$ is the elliptic bombon of type $(0,1)_2$ then, as we argued in general for the Case 2 in the previous section, there is a unique line through $x$ that touches $X$ only on $x$. Indeed, we may think of $X$ as the boundary of an ellipsoid $\mathcal E\subset\C^2$ and $T_xX$ is the closure of the unique affine complex line through $x$ within the (also unique) support real hyperplane of the ellipsoid $\mathcal E$ at $x$.

For the general case $n>2$, consider a non singular point $x$ in a bombon $X\subset\mathbb{CP}^n$. Since $T_xX\neq\mathbb{CP}^n$ is defined as a union of lines through $x$, to see that it is a hyperplane, it is enough to prove that for every plane $P$ through $x$, either  $P\cap T_xX$ is a line through $x$ or the whole of $P$.

By Lemma~\ref{sec}, $P\cap X$ is either a subspace or a bombon. If it is a subspace, it can be either $\{x\}$, a line through $x$ or $P$. In the three cases all the lines in $P$ through $x$ do not intersect $X$ in a circle, so that $P\subset T_xX$. If $P\cap X$ is a bombon, first observe that by the definition of the tangent space, $P\cap T_xX=T_x(P\cap X)$. By Theorem~\ref{teoDim2} we have three cases: $P\cap X$ is conical with $x$ its singular point or not, or else, $P\cap X$ is elliptic. In the first case, $T_x(P\cap X)=P$ and in the other two cases it is a line through $x$ by the case $n=2$.                \qed

\medskip

\noindent\emph{Proof of Theorem~\ref{thmEB}.} Let $X$ be an elliptic bombon in $\mathbb{CP}^n$. Consider $x\in X$, by Lemma~\ref{tangent}, $T_xX$ is a hyperplane that touches $X$ only in $x$ because for every line $L$ through $x$ in $T_xX$, we have that $L\cap X$ is not a circle nor $L$ (because $X$ is elliptic). Therefore, $X$ is the boundary of a complex ellipsoid in $\C^n\subset\mathbb{CP}^n$  by Proposition~\ref{propEll}.   \qed    

\smallskip
Observe that, according to the proofs of the previous section, we have proved more, because Case 2 was reduced to Case 1:

\begin{cor}\label{corell}
The bombon $X\subset \mathbb{CP}^n$ is an elliptic bombon if and only if there is a complex hyperplane $H$ that does not intersect $X$. Furthermore, any elliptic bombon in $\mathbb{CP}^n$ is algebraic of type $(0,n-1)_n$.
\end{cor}

\section{Hypersections of Bombons}

Recall that a section of a bombon is either a bombon or a subspace (Lemma~\ref{sec}). 
Our next  purpose  is to prove that elliptic and conical bombons  are precisely those bombons $X$ for  which the latter happens for \emph{hypersections}, that is, such that there is a hyperplane that intersects $X$ in a subspace.  

\begin{teo} \label{thmconical}
Let $X\subset \mathbb {CP}^n$ be a bombon. Suppose there is a hyperplane $H$ such that $H\cap X$ is a $k$-subspace, $-1\leq k\leq n-1$. Then $X$ is an elliptic bombon or a conical bombon. That is, $X$ is an algebraic bombon of  type $(0,m)_n$, where $0\leq m \leq n-1$.
\end{teo} 

\pf  Observe that the two extreme cases have been taken care of: for $k=-1$, by Corollary~\ref{corell}, and for $k=n-1$ by Theorem~\ref{thmdeg}. Now, we consider the general case where $X\subset\mathbb{CP}^n$ is a bombon with subspace hypersections.

Let $H$ be a hyperplane of $\mathbb{CP}^n$ such that $\Delta=H\cap X$ is a subspace of minimum dimension $k$ among all such possible subspace hypersections. We may further assume that $0\leq k<n-1$. 

Because $H\setminus \Delta$ is connected ($k<n-1$), then $H\setminus \Delta$ is contained in one of the two components of $\mathbb{CP}^n\setminus X$, say in $V$. Therefore, we have that $H\cap U=\emptyset$, where $U$ is the other component of $\mathbb{CP}^n\setminus X$. 

Now, we prove that the singular locus $\Sigma(X)$ is contained in $\Delta$. Indeed, if $x\in \Sigma(X)$, consider a line $L$ through $x$ that touches $U$. Because $x$ is singular, $L\setminus\{x\}\subset U$ and since any line touches a hyperplane in at least one point, then $L\cap H=\{x\}$ and $x\in \Delta$.

If $\Delta=\Sigma(X)$, let $\Gamma$ be a complementary subspace to $\Delta\neq\emptyset$ that passes through a point in $U$. Then, as in Lemma~\ref{SingJoinSmooth}, $\Gamma\cap X$ is a smooth bombon because $\Gamma$ intersects both $U$ and $V$.  And furthermore,  by Corollary~\ref{corell}, it is elliptic because $\Gamma\cap H$ is a hyperplane in $\Gamma$ that does not intersect $X$. So that $X$ is an algebraic bombon of type $(0,m)_n$ where $m=n-k-1$ is the dimension of $\Gamma$. We are left to prove that this is always the case, that is, that $\Delta=\Sigma(X)$.  

As in the Case 2 of Subsection~\ref{secEllBomI}, we identify $H$ with the hyperplane at infinity of $\mathbb{CP}^n$, and prove that  $U\subset \C^n$ is a convex set. Let  $a$ and $b$ be different points in $U$ and $L$ be the line through them. Then, $L\cap U$ es either an open halfplane or an open disk in $L\cap \C^n$, according to whether $L\cap H$ is in $\Delta$ or not, with $a$ and $b$ on the same side, $U$, and so is their segment. Therefore, $U$ and its closure $K=U\cup(X\setminus\Delta)\subset\C^n$ are convex sets. 

Consider a point $x\in X\setminus \Delta$, then, because $\Sigma(X)\subset\Delta$ and  Lemma~\ref{tangent}, $T_xX$ is a hyperplane. If $\Delta\not\subset T_xX$ then, using the convexity, we can push $T_xX$ parallel to itself and out of $K$ to a hyperplane $H^\prime$ for which $H^\prime\cap X=  \Delta\cap T_xX$, 
contradicting the minimality of the subspace hypersection.
So, for every $x\in X\setminus\Delta$, we have that $\Delta\subset T_xX$. Then, if $w\in\Delta$, all the lines through $w$ that touch another point of $X$ are contained in $X$ and so $w$ is a singular point. So $\Delta= \Sigma(X)$ and the proof is complete.					 \qed

\medskip

\begin{cor}
Any two  bombons in $\mathbb{CP}^{n}$ of type $(0,q)_n$, $0\leq q \leq n-1$  are algebraic bombons and therefore projetively equivalent. 
\end{cor}

\medskip
For further reference let us state the following corollary to the previous proof which resembles  Proposition~\ref{propEll}, but changes the hypothesis.

\begin{cor}\label{corEll} 
Let $X$ be a bombon in $\mathbb{CP}^n$. If there exists a hyperplane $H$ such that $X\cap H$ is a single non-singular point, then $X$ is elliptic, that is, it is algebraic of type $(0,n-1)_n$.   \qed
\end{cor}

\pf In the proof of Theorem~\ref{thmconical}, this case would yield that there exists a hyperplane that does not touch $X$ and therefore that $X$ is algebraic of type $(0,n-1)_n$.   \qed

\medskip
Lemma \ref{sec} and Theorem~\ref{thmconical} yield the following.

\begin{cor}\label{superraro}
Every hypersection of a non-elliptic smooth bombon is a bombon. \qed
\end{cor}

Now we prove the basic inductive step for the fullness of bombons theorem.

\begin{teo}\label{raro}
Let $X\subset \mathbb {CP}^n$ be a smooth bombon and let $H$ be a hyperplane of $\mathbb {CP}^n$.   Suppose that $H\cap X$ is a  full bombon of type $(p,q)_{n-1}$, with $p+q=n-3$, then $X$ is a full smooth bombon of type $(p+1,q+1)_n$.
\end{teo}

\pf Let $X\subset \mathbb {CP}^n$ be a smooth bombon with $U$ and $V$ the two components of $\mathbb {CP}^n\setminus X$. Suppose there is a complex hyperplane $H$ such that  $H\cap X$ is a full bombon of type $(p,q)_{n-1}$, with $p+q=n-3$. Then, $H\cap X=X^\prime\star\{x\}$, where $X^\prime$ is a bombon of type $(p,q)_{n-2}$.

Let $\mathcal{C}^\prime_U$ and $\mathcal{C}^\prime_V$ be cores of $X^\prime$. Then, $\Gamma_p=\mathcal{C}^\prime_U\star\{x\}$ and $\Gamma_q=\mathcal{C}^\prime_V\star\{x\}$ are subspaces of dimensions $p+1$ and $q+1$, respectivelly, that generate $H$ and intersect at $\{x\}$. 
Furthermore, $\Gamma_p \cap X=\{x\}$, $\Gamma_q \cap X=\{x\}$, $\Gamma_p \setminus \{x\}\subset U$ and $\Gamma_q \setminus \{x\}\subset V$.

Let $L$ be a line through $x$ not in $H$ and such that $L\cap X$ is a circle. The line  $L$ exists, otherwise $x$ would be a singular point of $X$ contradicting the fact that $X$ is smooth. 

Let $\Delta_p$ be the $(p+2)$-subspace generated by $\Gamma_p$ and $L$. 
Clearly, $X\cap \Delta_p$ is a bombon in $\Delta_p$ because $\Delta_p$ intersects both $U$ and $V$. 
Moreover, $\Gamma_p$ is a hyperplane of $\Delta_p$ with the property that $\Gamma_p\cap(X\cap \Delta_p)$ is the single point $x$ which is not singular because the line $L\subset\Gamma_p$ intersects $X$ in a circle. Consequently, by Corollary~\ref{corEll}, $X\cap \Delta_p$ is an elliptic bombon. 
Therefore, there is a $(p+1)$-subspace $\mathcal{C}_U$ contained in $U\cap \Delta_p\subset U$. 

If we interchange $p$ for $q$ and $U$ for $V$, verbatim, in the preceding paragraph, we obtain that there is a  $(q+1)$-subspace $\mathcal{C}_V$ contained in $V$. Therefore, because $(p+1)+(q+1)=n-1$, $X$ is a full smooth bombon of type $(p+1, q+1)_n$.
\qed

\section {All Bombons are Full}

\begin{teo}\label{unique}
If $X$ is a non elliptical smooth bombon, then for every $x\in X$, the tangent space $T_xX$ of the bombon $X$ at $x$, is a hyperplane and $\{x\}$ is the singular locus of  $T_x(X)\cap X$.
\end{teo}
\pf By Lemma~\ref{tangent}, the tangent space $T_xX$ is a hyperplane and by Corollary~\ref{superraro}, since $X$ is not elliptical, $T_X\cap X$ is a bombon and the point $x$ is in the singular locus of $T_xX\cap X$.  Suppose $\Sigma(T_xX\cap X)$ is not $\{x\}$, hence there is a line $L$ such that $x\in L\subset\Sigma(T_x(X)\cap X)$. Then, for every $y\in L$, we have that \begin{equation}\label{eq:Lyx}
 T_x(X)=T_y(T_x(X)\cap X)\subset T_y(X)\,,   
\end{equation} 
because $y$ is in the singular locus of $T_x(X)\cap X$. Consider $z\in(X\setminus T_x(X))$, and let $P$ be the plane generated by $L$ and $z$. We claim that $P \cap X$ is a bombon of type $(0,0)_2$. Indeed, it is not $P$ because the line from $z$ to $x\in L$ is not contained in $X$ ($z\not\in T_x(X)$), and it also contains a line: $L$. Therefore, $P \cap X$ is a bombon with a singular point $w\in L$ (because all the lines in a bombon of type $(0,0)_2$ pass through its singular point). Since $T_x(X)\subset T_w(X)$ by \eqref{eq:Lyx}, but $z\in T_w(X)\setminus T_x(X)$, then $T_w(X)=\mathbb{CP}^n$ and $w$ is a singular point of $X$ contradicting that it is smooth.
\qed

\smallskip

\begin{teo}\label{full}
Every bombon  $X\subset \mathbb{CP}^n$ is a full bombon. That is, $X=Y\star\Sigma(X)$, where  $\Sigma(X)$ is the singular locus of $X$ of dimension $k$ and  $Y$ is a smooth bombon of type $(p,q)_{n-k-1}$ with $p+q=n-k-2$.
\end{teo}
\pf The proof is by induction on the dimension $n$. If $n=2$, by Theorem~\ref{teoDim2} the theorem is true. Suppose it is true for $n-1$, we shall prove it for $n$. 
First of all, note that it is enough to prove the theorem when $X$ is smooth. Therefore, suppose $X$ is smooth. If $X$ is elliptical, then $X$ is full of type $(0,n-1)_n$. If $X$ is not elliptical, consider $x\in X$.  By Theorem~\ref{unique} $T_xX\cap X$ is a bombon and  $\{x\}$ is the only singular point of $T_x(X)\cap X$. Consequently, 
by induction $T_xX\cap X$ is full bombon of type $(p^\prime,q^\prime)_{n-1}$ with $p^\prime+q^\prime=n-3$. Then, by Theorem~\ref{raro} $X$ is a full bombon of type $(p,q)_n=(p^\prime+1,q^\prime+1)_n$. 
\qed

\smallskip 
An elliptical bombon in $\mathbb {CP}^n$ is a differential real $(2n-1)$-manifold embedded in a real $(2n)$-manifold.  Except por the singular points, the same is true for all bombons:

\begin{teo}\label{variedad}
Let $X\subset \mathbb {CP}^n$ be a bombon. Then $X\setminus \Sigma(X)$ is a real $(2n-1)$-manifold embedded  differentially in the real differentiable $(2n)$-manifold $\mathbb {CP}^n$. 
\end{teo}

\pf By Lemma~\ref{SingJoinSmooth} and the definition of join, it is clearly enough to prove the case of a smooth bombon. And observe that the Theorem is true for $n=2$ where both possible bombons (of type $(0,0)_2$ and $(0,1)_2$) are algebraic. 

Let $X\subset \mathbb{CP}^n$ be a smooth bombon and let $x\in X$ be any point.  Consider $\mathbb C^n\subset \mathbb {CP}^n$ as a standard differentiable chart ($\mathbb{CP}^n$ minus a hyperplane) such that $x\in \mathbb C^n=\R^{2n}$.  Then, $T_xX\cap \C^n$ is a real $(2n-2)$-subspace which is tangent to $X\cap \mathbb C^n$.  Let $L$ be a line through $x$ transversal to $T_xX$. By definition ($L\cap T_xX=\{x\}$), $L\cap X$ is a circle.  Let $\ell\subset L\cap \mathbb C^n$ be the real line tangent to 
$L\cap X \cap \mathbb C^n$ at $x$. Observe that $\ell\not\subset T_xX$, and define $\mathcal T_xX \subset \mathbb C^n$ to be the real $(2n-1)$-subspace generated by $T_xX\cap\mathbb C^n$ and the real line $\ell\subset \C^n$.   

We now show that the definition of $\mathcal T_xX$ is independent of the transversal $L$.   Let $L’$ be a different transversal 
line to $T_xX$ at $x$. Let $\ell’\subset L’\cap \mathbb C^n$ be the real line tangent to the circle $L’\cap X \cap \mathbb C^n$ at $x$.  We shall prove that 
$\ell’$ is a line of $\mathcal T_xX$.  To see this, let $P$ be the plane generated by $L$ and $L’$.  Then, $P \cap X$ is a bombon (in $\mathbb{CP}^2\cong P$) with $x$ not a singular point of $P \cap X$, because $L\subset P$.  This implies that $P \cap X\cap \mathbb C^n$ has a well defined unique tangent real $3$-dimensional plane $\Gamma$ at the point $x$. Then, both $\ell$ and $\ell’$ are in $\Gamma$ and furthermore, $\Gamma= P \cap \mathcal T_xX$. In particular,  $\ell’\subset \mathcal T_xX$.   

Next, we prove that any real line $\ell_1$ through $x$ in $\mathcal T_xX$ is  tangent to $X\cap \mathbb C^n$ at $x$.  Consider the complex line $L_1\subset \mathbb {CP}^n$ that contains $\ell_1$. If $L_1\subset T_xX$ there is nothing to prove. If not, $L_1\cap X$ is a circle and $\ell_1\subset (L_1\cap\C^n)$ is a real line tangent to it by the preceding paragraph.
\qed

\smallskip 

The following proposition includes an alternative proof of the manifold assertion.

\begin{prop}\label{propGrande}
Let $X$ be smooth bombon of type $(p,q)_n$,  and let $U$ and $V$ be the components of $\mathbb {CP}^n\setminus X$. Then, $p+q=n-1$ and $X$ is a real differentiable $(2n-1)$-manifold embedded in $\mathbb{CP}^{n}$. Moreover, 
\begin{enumerate}
\item Through every $u \in U$, there passes a core $p$-subspace  $C_u\subset U$ such that the inclusion $C_u\hookrightarrow U$ is a homotopy equivalence. Similarly, every $v \in V$, is in a core $q$-subspace  $C_v\subset V$ such that the inclusion $C_v\hookrightarrow V$ is a homotopy equivalence.
\item If $p,q>0$, then for every $x\in X$ there is a tangential $(n-1)$-subspace $H$ of $X$ such that $H\cap X$ is a bombon of class $(p-1,q-1)_{n-1}$, whose apex is the point $x$.
\end{enumerate}
\end{prop}

\pf First of all, by Theorem~\ref{full}, $p+q=n-1$ and we have that there exist subspaces $\mathcal{C}_U\subset U$, a core of $U$, and $\mathcal{C}_V\subset V$, a core of $V$, of respective dimensions $p$ and $q$.  

Let $\Gamma$ be any  $(p+1)$-subspace of $\mathbb {CP}^n$ containing $\mathcal C_U$.  Let $\{y\}=\Gamma\cap\mathcal C_V$ and observe that $y\in \mathcal C_V$ defines, and thus parametrizes all such $\Gamma$. We have that  $\Gamma\cap X$ is a bombon because $\Gamma$ intersects both $U$ and $V$. Then, by Corollary~\ref{corell}, $\Gamma\cap X$ is the boundary of a complex ellipsoid in $\Gamma$.
Since every $(p+1)$-dimensional subspace of $\mathbb {CP}^n$ containing $\mathcal C_U$ intersects $X$ in the boundary of a complex ellipsoid, then as in the proof of Lemma~\ref{cfb},  there is a fiber bundle with base space $\mathbb{CP}^q\cong\mathcal C_V$, fiber space $\mathbb S^{2p+1}$ and total space $X$ because any point $x\in X$ is in a unique $\Gamma=\mathcal C_U\star\{x\}$.  

This implies that $X$ is a real $(2n-1)$-dimensional manifold embedded in 
$\mathbb {CP}^n$.   Not only this, there is a fiber bundle with total space $V$, base space $\mathbb {CP}^q$ and fiber space $B^{2p+2}$. Furthermore, the subspace $\mathcal C_V$ intersects each fiber in a single point and therefore is a section of this bundle.  Hence, the inclusion $\mathcal C_V$ in $V$ is a homotopy equivalence.
Similarly, the inclusion $\mathcal C_U$ in $U$ is a homotopy equivalence.  

Now, we see that there are plenty of cores. Let $y\in U$ be such that $y\not\in\mathcal C_U$, and let $\Gamma$ be the $(p+1)$-subspace containing $\mathcal C_U$ and $y$. As above, $\Gamma \cap X$ is a an elliptic bombon, that is, $\Gamma \cap X$ is the boundary of a convex ellipsoid whose exterior is contained in $U$. Therefore, there is a $p$-subspace through $y$ contained in $U$. Similarly, for any $y\in V$, there is a $q$-subspace through $y$ contained in $V$.

Finally, if $p,q>0$ then $X$ is smooth but non elliptical. Then, by Theorem \ref{unique} 
for any $x\in X$, the tangent space $T_xX$ is  such that $T_xX\cap X$ is a bombon of type $(p-1,q-1)_{n-1}$ whose apex is the point $x$. 
\qed

\begin{cor}
Let $X$ be  bombon and suppose one of the two components of $\mathbb {CP}^n\setminus X$ is contractible. Then $X$ is an elliptical or conical bombon.
\end{cor}

\section{Homogeneity of smooth bombons}

We need some facts about the group of Moebius Transformations $\mathcal M$.

 Consider the circle $\hat{\R}$ which is the closure of $\R\subset\C\subset\mathbb{CP}^1$ and is naturally identified with the real projective line $\mathbb{RP}^1$. Its stabilizer, $St_\mathcal M(\hat{\R})$ (that is, the subgroup of $\mathcal M$ that fixes it as a set), is isomorphic to the projectivities of $\mathbb{RP}^1$. Its subgroup of index $2$ that preserves the orientation of $\hat{\R}$ acts as the hyperbolic motions of the Poincare upper half-plane model of the hyperbolic plane. Then, for every point $c$ in the upper half-plane, we have a well defined subgroup of hyperbolic rotations (parametrized by angles) with center $c$; and which, as Moebius transformations, also fix its conjugate $\bar{c}$. Therefore, given any circle $\Sigma$ in $\mathbb{CP}^1$ and a point $u\not\in\Sigma$, conjugating by an element in $\mathcal M$ that sends $\Sigma$ to $\hat{\R}$ and $u$ to the upper half-plane, there is a well-defined subgroup of \emph{rotations} of $\Sigma$ with \emph{center} $u$, denoted $\mathbb{S}_{\Sigma,u}=St_\mathcal M(\Sigma\cup \{u\})$ which is canonically isomorphic to $\mathbb{S}^1$; that is, there is a well defined isomorphism $$\phi_{\Sigma,u}:\mathbb{S}^1\to\mathbb{S}_{\Sigma,u}\,.$$ 
 These rotations also fix the \emph{conjugate} of $u$ with respect to $\Sigma$, $v$ say, in the other component (than $u$) of $\mathbb{CP}^1\setminus\Sigma$, and the two isomorphisms from $\mathbb{S}^1$ to $\mathbb{S}_{\Sigma,u}=\mathbb{S}_{\Sigma,v}\subset\mathcal M$ are inverse to each other (the only non trivial automorphism of $\mathbb{S}^1$).

\begin{lema}\label{lema:acciones}
Let $X$ be a smooth bombon with components $U$ and $V$. Given cores of $X$, $\mathcal C_U\subset U$ and $\mathcal C_V\subset V$, there is a $\mathbb{S}^1$-action on $X$ that extends to $U$ and makes $X$ the total space of a principal $\mathbb{S}^1$-bundle over $\mathcal C_U\times \mathcal C_V$.
\end{lema}

\pf
Given $u\in \mathcal C_U$ and $v\in \mathcal C_V$, the line through them, $L=u\star v$, intersects $X$ in a circle, $\Sigma$, because it intersects the two components of $\mathbb{CP}^n\setminus X$. Then, $\Sigma\subset X$ has a $\mathbb{S}^1$-action as rotations along $\Sigma$ with center $u$, that extends to the disk in $L$ with boundary $\Sigma$ and containing $u$.
Since every $x\in X$ is in a unique such circle (indeed, take $u=(x\star \mathcal C_V)\cap \mathcal C_U$ and $v=(x\star \mathcal C_U)\cap \mathcal C_V$), this completes the proof, because, likewise, every point in $U\setminus \mathcal C_U$ is in a unique such line, and the points in $\mathcal C_U$ are fixed points of the $\mathbb{S}^1$-actions on lines that have been considered.
\qed

\begin{teo}\label{teo:homogeneidad}
A smooth bombon is homogeneous.
\end{teo}

\pf
Consider two points $x$ and $y$ in a smooth bombon $X$ with componenets $U$ and $V$. If the line through them, $L=x\star y$,  intersects $X$ in a circle, consider cores $\mathcal C_U\subset U$ and $\mathcal C_V\subset V$ that touch $L$ (they exist by Proposition~\ref{propGrande}(1)). Then, by the proof of Lemma~\ref{lema:acciones}, $x$ and $y$ are in the same orbit of the $\mathbb{S}^1$-action on $X$ defined by $\mathcal C_U$ and $\mathcal C_V$, and therefore there is a homeomorphism of $X$ (a rotation, one might say) that sends $x$ to $y$.

If the line $L$ through $x$ and $y$  is contained in $X$, consider $z\in X$  such that the lines from $z$ to $x$ and to $y$ are not in $X$; it exists because $X$ is smooth so that $X\not\subset T_x(X)\cup T_y(X)\neq \mathbb{CP}^n$. Therefore, there are homeomorphisms of $X$ that send $x$ to $z$ and $z$ to $y$, so that their composition sends $x$ to $y$.  
\qed

\section{Apendix: The proof of Theorem \ref{thmB}}

Theorem \ref{thmB} is an unexpected  characterization of complex ellipsoids with no analogue over the real numbers. Some of the ideas in \cite{ABMCE} for dimension $2$, work nicely for all $n\geq 2$, and that is what we develop here. 

First note that all $1$-dimensional complex ellipsoids are disks. Therefore, Theorem \ref{thmB} is true for $n=1$.
 
An \emph{abstract linear space} consists of a set $\Lambda$  together with a distinguished family of subsets, called \emph{abstract lines}, satisfying the following property:  given different $x,y\in \Lambda$, there is a unique abstract line containing $x$ and  $y$.  Typical examples of abstract linear spaces are euclidean $n$-space, projective spaces and complex $n$-spaces. 

Our interest lies in the abstract linear space $LS^{2n-1}$, where $\Lambda=\mathbb S^{2n-1}$, the unit sphere of $\mathbb C^n$, and an abstract line is the intersection of a complex line of  $\mathbb C^n$ with $\mathbb S^{2n-1}$, which is a circle in $\mathbb S^{2n-1}$. 

A subset $A$ of an abstract linear space $\Lambda$ is \emph{linearly closed} if for any $x,y\in A$, the abstract line through $x$ and $y$ is contained in $A$.  Since the intersection of linearly closed subsets is linearly closed, given any $Y\subset \Lambda$, there exists a unique minimal linearly closed subset containing $Y$, called its \emph{linear closure}.

\begin{lema}\label{lemclo}
Suppose $K\subset\mathbb{C}^n$ is a convex set such that every complex line intersects $K$ either in the empty set, a single point or a disk, and  with the property that the complex ellipsoid of minimal volume containing it is the unit ball. Then $K\cap \mathbb S^{2n-1}$ is a linearly closed subset of the abstract linear space $LS^{2n-1}$.
\end{lema}
\pf 
Let $x,y\in K\cap \mathbb S^{2n-1}$, $x\neq y$, and let $L$ be the complex line through $x$ and $y$. We have to prove that $L\cap \mathbb S^{2n-1}\subset K$. 
Let $B^{2n}$ be the unit ball in $\mathbb{C}^n$, whose boundary is $\mathbb S^{2n-1}$. By hypothesis, $L\cap K$ and $L\cap B$ are two disks in $L$ which share two different points, $x$ and $y$ on their boundary. Then,  
$ K \subset B^{2n}$ implies that $L\cap K = L\cap B^{2n}$, and therefore that $L\cap \mathbb S^{2n-1}\subset K$. \qed

\smallskip

\begin{prop}\label{prop:lin_cl}
The linearly closed proper subsets of the abstract linear space $LS^{2n-1}$ are the intersections of $\mathbb S^{2n-1}$ with complex affine $k$-subspaces. 
\end{prop}

By induction, it will be enough to prove the following:  \emph{let $A$ be a linearly  closed subset of the abstract linear space $LS^{2n-1}$ and suppose $H$ is a complex hyperplane intersecting the interior of $B^{2n}$.  If $H\cap \mathbb S^{2n-1} \subset A$ but there is $a\in A\setminus (H\cap \mathbb S^{2n-1})$, then 
$A=LS^{2n-1}$.}

The proof of the above claim follows immediately from the next two lemmas. 
For which we need to make precise some standard definitions. 
A topological space $U$ is topologically embedded in $\mathbb C^n$, if there is a continuous injective map $f: U\to\mathbb C^n$. A topological space $A$ is locally embedded in $\mathbb C^n$ if for every point $x\in A$, there is a neighborhood $U$ of $x$ in $A$ which can be topologically embedded in $\mathbb C^n$.

\begin{lema}
If $A\subset\mathbb S^{2n-1}$  is closed and $A\neq\mathbb S^{2n-1}$, then $A$ can be embedded in $\mathbb C^{n-1}$.
\end{lema}

\pf Since $A$ is proper, there exists a point $w \in\mathbb S^{2n-1}$ such that $w\notin A$.  Let $L$ be the complex line through $w$ and the origin and let  $L^\perp$ be the orthogonal complex hyperplane to $L$ at the origin, which we identify with $\mathbb C^{n-1}$. 

Let $\pi:\mathbb S^{2n-1}\setminus \{w\} \to L^\perp$  be the geometric projection from $w$, that is, given $x\in \mathbb S^{2n-1}$ different  from $w$, $\pi(x)$ is the intersection with $L^\perp$ of the complex line through $x$ and $w$; well defined because this complex line is not parallel to  $L^\perp$.

We claim that $\pi |_A$ is an embedding. Indeed, if $x, y \in A$ are such that $\pi(x)=\pi(y)$ then their abstract lines through $w$ coincide. If $x$ and $y$  were different points, it would imply that $w\in A$ because $A$ is linearly closed, which contradicts the choice of $w$.        \qed 

\medskip
From now on let $A$ be a linearly  closed subset of the abstract linear space $LS^{2n-1}$ and suppose $H$ is a complex hyperplane intersecting the interior of $B^{2n}$ such that $H\cap \mathbb S^{2n-1} \subset A$ but there is $a\in A\setminus (H\cap \mathbb S^{2n-1})$. 

\begin{lema}
The set  $A$ is not locally embedded in $\mathbb C^{n-1}$.
\end{lema}

\pf By hypothesis $A$ contains $H\cap \mathbb S^{2n-1}$, where $H$ is a complex hyperplane intersecting the interior of $B^{2n}$.
Furthermore,  there is yet another point $a\in A\setminus H$.

For every $x\in H\cap \mathbb S^{2n-1}$, let $L_x$ be the complex line through $a$ and $x$.  Note that $L_x\cap H=\{x\}$, and if $x\not=x'\in H\cap \mathbb S^{n-1}$,
then $L_x\cap L_{x'}=\{a\}$, because two complex lines intersect in at most one point. Let 
$$\Omega=\bigcup_{x\in H\cap \mathbb S^{2n-1}}L_x \subset A\subset \mathbb S^{2n-1}.$$  
Since every abstract line of $LS^{2n-1}$ is actually a circle, then 
$$\Omega\setminus\{a\}=\bigcup_{x\in H\cap \mathbb S^{2n-1}}(L_x\setminus\{a\})$$ 
is an $\R$-bundle over $\mathbb S^{2n-3}$.  

If $n\geq 2$, this bundle must be trivial and therefore $\Omega$ is homeomorphic to a closed cylinder  $S^{2n-3}\times [0,1]$ modulo its boundary,
$\frac{\mathbb S^{2n-3}\times [0,1]}{\mathbb S^{2n-3}\times\{0,1\}}$, where the point $a\in \Omega$ is precisely $\frac{\mathbb S^{2n-3}\times \{0,1\}}{\mathbb S^{2n-3}\times\{0,1\}}$.
  Consequently, any neighborhood of $a$ in $\Omega$ cannot be topologically embedded in $\mathbb C^{n-1}$, because it contains the cone of two disjoint $(2n-3)$-spheres. Since $\Omega \subset A$, then $A$ is not locally embedded in $\mathbb C^{n-1}$.  
 If  $n=2$ and this  $\R$-bundle over $\mathbb S^1$ is trivial, the same argument shows that $A$ is not locally embedded in $\mathbb C$. This time, the $\R$-bundle over $\mathbb S^1$ may be not trivial, but hence $\Omega$ is homeomorphic to the projective plane 
 This is imposible because $\Omega\subset A\subset \mathbb S^3$.
  \qed

\bigskip

\noindent \emph{Proof of Theorem~\ref{thmB}.} 
 Suppose $K\subset\mathbb{C}^n$ is a convex set such that every complex line intersects $K$ either in the empty set, a single point or a disk and let 
$E$ be the complex ellipsoid of minimal volume containing $K$ (see, e.g., \cite{ABM1}). We may assume without loss of generality that $E$ is the unit ball of $\mathbb{C}^n$. We must prove that $K = E$, which is equivalent to $K\cap \mbox{bd} E=K\cap \mathbb S^{2n-1}=\mathbb S^{2n-1}$.

By Lemma \ref{lemclo},  $K\cap \mathbb S^{2n-1}$ is a linearly closed subset of $LS^{2n-1}$ and by Proposition\ref{prop:lin_cl}, we have to prove that $K\cap \mathbb S^{2n-1}$ is not the intersection of $\mathbb S^{2n-1}$ with a complex affine $k$-subspace $H$. This is so because by John's Theorem \cite{John} (see also \cite{JMM}),
if $K$ is a convex set which has the unit ball $B^{2n}$ as its minimal complex ellipsoid, then $K\cap \mathbb S^{2n-1}$ cannot be contained in a proper complex hyperplane of $\mathbb C^n$.  \qed

\medskip
As in the case of elliptic bombons, we think that all bombons are algebraic.
\bigskip

\noindent {\bf Acknowledgments.} Luis Montejano acknowledges  support  from  PAPIIT-UNAM with project 35-IN112124. Javier Bracho acknowledges  support  from projects PAPIIT-UNAM  IN109023 and SECIHTI  CBF-2025-I-224.

\end{document}